\documentclass[11pt]{article}
\usepackage{amsfonts,amsmath,amssymb,amsthm, amscd}
\usepackage[dvips]{epsfig,graphics}

\numberwithin{equation}{section}

\newtheorem{theorem}{Theorem}[section]

\newtheorem{lemma}[theorem]{Lemma}

\newtheorem*{theoremB}{Theorem of Burde and de Rham}

\theoremstyle{definition}
\newtheorem{definition}[theorem]{Definition}
\newtheorem{example}[theorem]{Example}
\newtheorem{remark}[theorem]{Remark}

\newcommand{\D}{\Delta}
\DeclareMathOperator{\im}{Im}

\def\<{{\langle}}
\def\>{{\rangle}}

\def\a{{\alpha}}

\def\d{{\delta}}
\def\g{{\gamma}}
\def\Z{\mathbb Z}

\def\C{\mathbb C}
\def\a{\alpha}

\def\L{{\Lambda}}

\def\e{\epsilon}
\def\A{{\cal A}}

\def\o{\otimes}

\begin{document}

\title{On a Theorem of Burde and de Rham}

\author{Daniel S. Silver \and Susan G. Williams\thanks{Both authors partially supported by NSF grant
DMS-0706798.} \\ {\em
{\small Department of Mathematics and Statistics, University of South Alabama}}}

\maketitle 

\begin{abstract}  We generalize a theorem of Burde and de Rham characterizing the zeros of the Alexander polynomial. Given a representation of a knot group $\pi$, we define an extension $\tilde \pi$ of $\pi$, the Crowell group. For any ${\rm GL}_N \C$ representation of $\pi$, the zeros of the associated twisted Alexander polynomial correspond to representations of $\tilde \pi$ into the group of dilations of $\C^N$. 

\end{abstract} 

A classic theorem of G. Burde, and independently G. de Rham
characterizes the zeros of the Alexander polynomial of a knot in terms of representations of the knot group by dilations, certain affine transformations of the complex plane. More recently, twisted Alexander polynomials, which incorporate information from linear representations of the knot group, have proved useful in many areas (see the survey \cite{friedl}). We generalize the theorem of Burde and de Rham for twisted Alexander polynomials, replacing  dilations of the complex plane with dilations of $\C^N$, where $N$ is the dimension of the representation. The knot group is replaced by a natural extension, 
the Crowell group,  introduced in \cite{crowell}.

\section{Theorem of Burde and deRham}

\begin{definition} \label{dil.def} A {\it dilation} of $\C^N$ is a transformation $d : \C^N \to \C^N$ of the form $d(z) = \a z+ b$, for some $\a \in \C^*= \C\setminus \{0\}$ and $b \in \C^N$. 

\end{definition}

The number $\a$ is called the {\it dilation ratio} of $d$. Dilations of $\C^N$ form a group ${\cal D}_N$ under composition. It is easy to see that conjugate elements have the same dilation ratio $\a$. 

Assume that $k \subset {\mathbb S}^3$ is a knot with group $\pi = \pi_1({\mathbb S}^3 \setminus k)$ given by a Wirtinger presentation 
$\<x_0, x_1, \ldots, x_n \mid r_1, \ldots, r_n\>$.  Let $\e: \pi \to \<t \>\cong \Z$ be the abelianization homomorphism sending each  $x_i$ to $t$. Wirtinger relations have the form $x_i x_j = x_k x_i$. Any one relation is a consequence of the others, and so we omit a single relation.

The {\it Alexander module} is the $\Z[t^{\pm 1}]$-module $\A$ with generators $x_0, x_1, \dots, x_n$ and $n$ relations of the form 
\begin{equation} \label{alex.relation} x_i + t x_j = x_k + t x_i \end{equation} corresponding to the Wirtinger relations. The {\it based Alexander module } $\A^0$ is the quotient $\A/\<x_0\>.$ Both $\A$ and $\A^0$ are invariants of the knot. 

The module $\A^0$ has an $n \times n$ presentation matrix $M(t)$. 
The greatest common divisor of the $(n-r+1) \times (n-r+1)$ minors is $\D_r(t)$, the $r$th {\it Alexander polynomial} of $k$.  It is well defined up to multiplication by a unit in $\Z[t^{\pm 1}]$; we  normalize so that it is a polynomial with
$\D_r(0) \ne 0$. We denote $\D_1(t)$ more simply by $\D(t)$, and refer to it as the {\it Alexander polynomial} of $k$.
It is well known that for each $r$, the set of zeros of $\D_r(t)$ is closed under inversion.

The set of representations $\d: \pi \to {\cal D}_1$ is denoted by  ${\rm Hom}(\pi, {\cal D}_1)$. Such  representations have the form $\d(x_i): z \mapsto \a z +b_i$. 
For fixed $\a\in \C^*$, we denote by ${\rm Hom}(\pi, {\cal D}_1)_\a$ the subset of ${\rm Hom}(\pi, {\cal D}_1)$ consisting of representations $\d$ such that $\d(x_i)$ has dilation factor $\a$, for each $i$. Any such $\d$ is uniquely determined  by the vector $(b_0, \ldots, b_n) \in \C^{n+1}$.

Representations $\d\in {\rm Hom}(\pi, {\cal D}_1)_\a$ with cyclic images correspond to vectors of the form $(b_0, b_0, \ldots, b_0)$. In order to remove them from consideration, we restrict our attention to the subset ${\rm Hom}(\pi, {\cal D}_1)_\a^0$ consisting of {\it based representations}, those for which $b_0=0$. For $\a\ne 1$, every representation in ${\rm Hom}(\pi, {\cal D}_1)_\a$ is conjugate in ${\cal D}_1$ to a based one, via the translation $z \mapsto b_0/(1-\a)$.  It is easily seen that the set of vectors $(b_1, \ldots, b_n) \in \C^{n}$ arising from elements of ${\rm Hom}(\pi, {\cal D}_1)_\a^0$  is the nullspace of the matrix $M(\a^{-1})$.  We identify ${\rm Hom}(\pi, {\cal D}_1)_\a^0$ with this subspace.  The first part of the following theorem follows from well-known fact that for each $r$, the set of zeros of $\D_r(t)$ is closed under inversion.

\begin{theoremB} [\cite{burde}, \cite{rham}] \label{br.theorem} \emph{(1)} For any $\a \in \C^*$, the dimension of ${\rm Hom}(\pi, {\cal D}_1)_\a^0$ is equal to $\max \{r \mid \D_r(\a) =0 \}$. 

\emph{(2)} When $\a$ is a zero of an irreducible factor of $\D(t)$ that does not occur in $\D_2(t)$, any two representations $\d\in {\rm Hom}(\pi, {\cal D}_1)_\a^0$ are conjugate in ${\cal D}_1$.
\end{theoremB}

\section{Twisted Alexander invariants} Let $\g: \pi \to {\rm GL}_N \C$ be a representation. Denote the image $\g(x_i)$ by $X_i$, for $0 \le i \le n$.
We obtain a representation $\e\o \g: \pi \to {\rm GL}_N \L$, where $\L$ is the Laurent polynomial ring $\C[t^{\pm 1}]$,  by mapping each $x_i$ to $t X_i$. Let $M_N\L$ denote the ring of all $N\times N$ matrices over $\L$. The $\g$-{\it twisted Alexander module} $\A_\g$ is the $M_N\L$-module with generators
$x_0, \ldots, x_n$ and relations

\begin{equation}\label{twisted.alex.relation} x_i + (t X_i) x_j = x_k + (t X_k) x_i \end{equation} corresponding to the Wirtinger relations for $\pi$. As in the untwisted case, any relation is a consequence of the remaining ones. Hence we omit a single relation. 

The {\it based $\g$-twisted Alexander module} $\A_\g^0$ is the quotient $\A_\g/\<x_0\>$. Both $\A_\g$ and $\A_\g^0$ are invariants of $k$. (Homological interpretations of $\A_\g$ and $\A_\g^0$ are given in \cite{swtwisted}.)

We regard $M_\g$ as 
an $nN \times nN$-matrix with entries in $\L$. The greatest common divisor of the $(nN-r+1) \times (nN-r+1)$ minors, denoted by $D_{\g, r}(t)$, is well defined up to multiplication by a unit in $\L$. We call $D_{\g, 1}(t)$ the {\it Alexander-Lin polynomial} (see {\cite{sw09}), and denote it more simply by $D_\g(t)$. 

\begin{remark} $D_\g(t)/\det (tX_0 -I)$ is the invariant $W_\g(t)$ defined by Wada  \cite{wada}. In general, it is a rational expression. 

Another version of the twisted Alexander polynomial, defined homologically, was given by P. Kirk and C. Livingston in \cite{kl}. Their invariant is $D_\g(t)/f(t)$, where $f(t)$ is a certain factor of $\det (tX_0 -I)$ (see \cite{kl} for details). 
\end{remark}

\section{The Crowell group of a knot} In \cite{crowell}, R.H. Crowell introduced a group extension of a knot group $\pi$ with the useful property that it abelianizes to the Alexander module. A twisted generalization was presented in \cite{swderived}.

\begin{definition} \label{derived} \cite{crowell} The  {\it derived group} $\tilde G$ of an action $S \times G \to G$ has generator set $S \times G$, written  $\{^sg \mid s \in S, g\in G\}$, and relations
$^s(gh) = (^s g)( ^{sg}h)$, for all $g, h \in G, s \in S$. \end{definition}

\begin{remark} (1) Crowell denotes the element $^sg$ by $s\wedge g$. We prefer the exponential notation because it is more compact and it suggests the operator notation
that we employ below. 

In \cite{swderived}  the  element $^sg$ is written as $g^s$. Here the mild change of notation is intended to avoid possible sources of notational confusion.

(2) As Crowell remarked \cite{crowell}, the derived group $\tilde G$ of an action should not be confused with commutator subgroup $[G, G]$. The term ``derived" 
refers to free group derivations. \end{remark}

\begin{lemma} \label{derivedlemma} \cite{crowell} For the derived group $\tilde G$ of any action $S \times G \to S$, the following hold.
 \item\emph{(1)} $^s\!g = 1$ if and only if $g =1$;
\item\emph{(2)} $({}^s\!g)^{-1} = {}^{sg}(g^{-1})$;
\item\emph{(3)} every nontrivial element of $\tilde G$ has a unique expression of the form $^{s_1}g_1 \cdots {}^{s_k}g_k$, where each $g_i \ne 1$ and $s_{i+1} \ne s_i g_i$ for $i=1, \ldots, k-1$.

\end{lemma}
   
The complete proof of Lemma \ref{derivedlemma} can be found in \cite{crowell}. The main idea behind the proof of {\it (3)} is that any product of the form ${}^{s_i}g_i\ {}^{s_ig_i}g_{i+1}$ can be rewritten as $^{s_i}(g_ig_{i+1})$.

\begin {example} \label{ex} 


Consider a knot group with Wiritinger presentation $\pi= \<x_0, x_1, \ldots, x_n \mid r_1, \ldots, r_n\>$. 
We describe the derived group $\tilde \pi$ for two types of action.

(1) Consider the abelianization $\pi \to \<t\>$ mapping each Wirtinger generator to $t$. The group $\pi$ acts on 
$S=\<t\>$ via $t^\nu x_i = t^{\nu+1}$, for each $\nu \in \Z$ and $0\le i \le n$. The derived group $\tilde \pi$ has generators 
$^{t^\nu}\!x_i$, $1\le i \le n,\ \nu \in \Z$. Relations have the form 
$${}^{t^\nu}\!\!x_i \ {}^{t^{\nu+1}} \!\!x_j = {}^{t^\nu}\!\!x_k\ {}^{t^{\nu+1}} \!\!x_i.$$
Writing $^{t^\nu}\!x_i$ as $x_{i, \nu}$, a
notation that arises in the Reidemeister-Schreier process (see for example \cite{rap} or \cite{mks}), the relations take the form: 
$$x_{i, \nu} x_{j,  \nu+1} = x_{k, \nu} x_{i, \nu+1}.$$

The abelianization of $\tilde \pi$ has the structure of a $\Z[t^{\pm 1}]$-module, the action of $t$ corresponding to conjugation in $\pi$ by $x_0$. It is not difficult to see that the module is isomorphic to the Alexander module
$H_1(\tilde X, \tilde b_0)$, where  $\tilde X \to X$ is the infinite cyclic cover of the knot exterior, $\tilde b_0$ is the preimage of a base point $b_0$. 
This example was the main motivation behind Crowell's construction. 

(2) Assume that $\g: \pi \to  {\rm GL}_N \C$ is a linear representation. Let $S$ be the image of $\g$. Then $\pi$ acts on $S$ via $sg = s\g(g)$, for each $g \in \pi$ and $ s\in S$. In the derived group $\tilde \pi$, each Wirtinger generator $x_i$ becomes a family $^Sx_i$ of generators. Similarly, each Wirtinger relation $r_i$ becomes a family $^Sr_i$ of relations. Any family of relations is a consequence of the remaining ones, and so we omit a single family. We have:

\begin{equation}\label{crowell.pres} \tilde \pi = \<{}^S\!x_0, \ldots, {}^S\!x_n \mid {}^S\!r_1, \ldots, {}^S\!r_n\> \end{equation}

More explicitly, any Wirtinger relation $x_i x_j = x_k x_i$ of $\pi$ gives the relation $^s(x_i x_j) = {}^s(x_k x_i)$ in $\tilde \pi$, for each $s \in S$. Writing $X_i$ and  $X_k$ for $\g(x_i)$ and $\g(x_k)$, respectively, these relations can be rewritten as 
$^sx_i\ ^{sX_i}x_j = \!^s x_k\ ^{sX_k}x_i$.

In the case that $\g$ is trivial (so that $S$ is the singleton), $\tilde \pi$ is isomorphic to $\pi$.

\end{example}

\begin{definition} If $\g: \pi \to  {\rm GL}_N \C$ is a representation of the group of a knot $k$, then the {\it Crowell group of $k$ (with respect to $\g$)} is the derived group  $\tilde \pi$ of the associated group action. \end{definition}

Example \ref{ex} (1) suggests Crowell's topological interpretation of the  derived group $\tilde G$ (cf. \cite{crowell}). Let $B, *$ be a connected space and base point such that $\pi_1(B, *) \cong G$. Let $p: E \to B$ be the covering space corresponding to the permutation representation $S \times G \to S$. The preimage $p^{-1}(*)$ is identified with $S$. Following \cite{crowell}, we denote by 
$s \wedge g$ the relative homotopy class of paths lying above $g$ beginning at $s$  and ending at $s g$. The set of such homotopy classes is a groupoid under concatenation, with
$$ (s \wedge  g)(sg \wedge h) = s \wedge g h,$$
for $g, h \in G$ and $s \in S$.  We define an analogous groupoid structure on $S \times G$ by 
$$ (s, g)(sg, h) = (s, g h).$$
The derived group $\tilde G$ is the smallest group containing this groupoid. More preciesely, the function $\psi: S \times G \to \tilde G$ mapping $(s,g) \mapsto {}^sg$ is a groupoid morphism  with the universal property that for any group $A$ and groupoid morphism $f: S \times G \to A$, there is a unique homomorphism $\bar f: \tilde G \to A$ such that $\bar f \psi = f$. 

\begin{figure}
\begin{center}
\includegraphics[height=2.5 in]{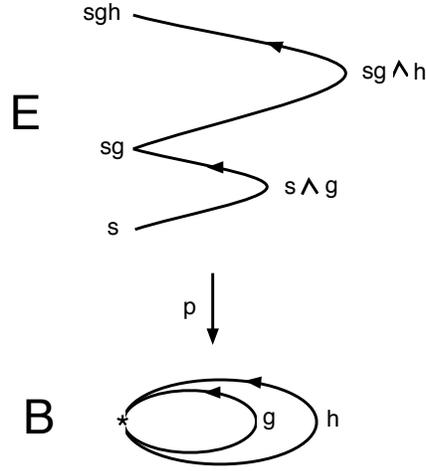}
\caption{$(s \wedge  g)(sg \wedge h) = s \wedge g h$}
\label{cover}
\end{center}
\end{figure}

Example \ref{ex} (2) reveals additional structure in a derived group whenever $S$ is a monoid and the action satisfies an additional axiom
\begin{equation}\label{axiom} (st)\cdot g = s(t\cdot g),\ {\rm for\ all}\ s, t \in S, g \in G.\end{equation}
  In this case, $\tilde G$ is an $S$-group, a group with operators. The general notion was developed by E. Noether \cite{noether}, who attributed the idea to Krull.  In general, an $S$-{\it group} is a group $\tilde G$ together with a set $S$ ({\it operator set}) and a function $S \times \tilde G \to \tilde G$, written $(s, \tilde g) \mapsto s \tilde g$, such that for each fixed $s$, the restricted map $\tilde g \mapsto s  \tilde g$ is an endomorphism of $\tilde G$. When the operator set $S$ is empty, an $S$-group is a group in the usual sense.

\begin{lemma}\label{sgroup} The derived group $\tilde G$ of an action of $G$ on a monoid $S$ satisfying condition \emph{(\ref{axiom})} is an $S$-group. \end{lemma} 

\begin{proof} For any $s \in S$ and $\tilde g={}^{s_1}\!g_1 \cdots {}^{s_k}\!g_k \in \tilde G$,  we define $s ({}^{s_1}\!g_1 \cdots {}^{s_k}g_k)$ to be ${}^{ss_1}\!g_1 \cdots {}^{ss_k}\!g_k$. Since
$$s\,{}^t\!(gh) = {}^{st}\!(gh) = ({}^{st}\!g) ({}^{(st)g}h) = ({}^{st}\!g) ({}^{s(tg)}\!h) = s ({}^t\!g\, {}^{tg}\!h)$$
for all $s, t \in S,\ g, h \in G$, the action is well defined. 

For fixed $s \in S$, the mapping $\tilde g \mapsto s \tilde g$ is  obviously a homomorphism. 
\end{proof}

Many notions from group theory extend to $S$-groups in a natural way.
Much like modules, $S$-groups can often be described in a simple fashion using $S$-group generators $x_i$ and relators $r_j$. For this, we regard $x_i$ as ${}^e\! x_i$, where $e\in S$ is the identity element. In a presentation, $x_i$ represents the orbit ${}^S\! x_i$. 
Similar considerations apply to $r_j$. For the $S$-group structure described in Lemma \ref{sgroup}, any group presentation of $G$ conveniently serves as an $S$-group presentation of $\tilde G$.

An $S$-{\it group homomorphism} (or {\it $S$-group representation}) $\d: G \to H$ is a mapping between $S$-groups that is a homomorphism of the underlying groups respecting the $S$-action in the sense that $\d(s g) = s \d(g)$, for all $g \in G, \ s \in S$. 


\section{Twisting the Theorem of Burde and de Rham} Let $\g: \pi \to  {\rm GL}_N \C$ be a representation of the knot group. Recall from Definition \ref{dil.def} that elements of ${\cal D}_N$ have the form 
$d(z) = \a z+ b$, for for some $\a\in \C^*$ and $b \in \C^N$. There is a left action of ${\rm GL}_N \C$ on ${\cal D}_N$: \begin{equation} \label{action1} (A, d) \mapsto Ad, \end{equation} where \begin{equation} \label{action2} Ad: z \mapsto \a z + Ab. \end{equation} By restriction, 
$S= \im\g$ acts as well, and it is easy to check that ${\cal D}_N$ is an $S$-group. 

Consider an $S$-group representation $\d: \pi_\g \to {\cal D}_N$. As in the untwisted case, any two $S$-group generators $x_i, x_j$ of $\pi_\g$ are conjugate, and hence their images under $\d$ must have the same dilation ratio $\a$.

Motivated by the untwisted case, we say that an $S$-group representation is {\it based} if $b_0 \in \C^N$ is the zero vector. (Note, however, that conjugation by a translation does not preserve the set of $S$-group representations in general.) We denote by ${\rm Hom}(\pi_\g, {\cal D}_N)^0_\a$ the set of based $S$-group representations $\d: \pi_\g \to {\cal D}_N$ such that each $\d(x_i)$ has dilation factor $\a$. Any $ \d \in {\rm Hom}(\pi_\g, {\cal D}_N)^0_\a$ is uniquely determined by the vector $(b_1, \ldots, b_n) \in (\C^N)^n \cong \C^{nN}$ with $\d(x_i): z\mapsto \a z +b_i$, for each $i$.

When does a vector $(b_1, \ldots, b_n)  \in \C^{nN}$ correspond to a representation in ${\rm Hom}(\pi_\g, {\cal D}_N)^0_\a$?  Each relation $x_i \ {}^{X_i} x_j= x_k\  {}^{X_k}x_i$ requires that $$\d(x_i)\  ({X_i} \d) (x_j)= \d(x_k)\ ({X_k}\d)(x_i).$$
From the form (\ref{action2}), the condition is equivalent to 
\begin{equation}\label{equivalent} \a b_i +  X_i b_j = \a b_k +  X_k b_i \end{equation} Other relations ${}^Ax_i\ {}^{AX_i}x_j = {}^Ax_k\  {}^{AX_k}x_i$, where $A\in S$ is arbitrary, are immediate consequences. 

The linearity of  condition (\ref{equivalent}) implies that the set of vectors $(b_1, \ldots, b_n)$ corresponding to representations in $ {\rm Hom}(\pi_\g, {\cal D}_N)^0_\a$ is a subspace of $\C^{nN}$. We identify $ {\rm Hom}(\pi_\g, {\cal D}_N)^0_\a$ with the subspace.

\begin{theorem} \emph{(1)} For any $\a \in \C^*$, the dimension of ${\rm Hom}(\pi_\g, {\cal D}_N)_\a^0$ is equal to $\max \{r \mid \D_{\g, r}(\a^{-1}) =0 \}$. 

\emph{(2)} When $\a$ is a zero of an irreducible factor of $\D_\g(t)$ that does not occur in $\D_{\g, 2}(t)$, any two nontrivial representations of 
${\rm Hom}(\pi_\g, {\cal D}_N)_\a^0$ are conjugate in ${\cal D}_N$.
\end{theorem}

\begin{proof} Assume that $(b_1, \ldots, b_n)$ determines a representation in ${\rm Hom}(\pi_\g, {\cal D}_N)_\a^0$. Since $\a \ne 0$, we can rewrite condition (\ref{equivalent}): 
\begin{equation}\label{equiv1} b_i + ( \a^{-1}X_i) b_j = b_k +  (\a^{-1}X_k)b_i . \end{equation} Comparing this condition with (\ref{twisted.alex.relation}), we see that $(b_1, \ldots, b_n)$ is in the null space of the relation matrix $M_\g$ for $\A_\g^0$ evaluated at $t=\a^{-1}$. Hence $D_\g(\a^{-1})=0$. In this case, the possible vectors $(b_1, \ldots, b_n)$ form a subspace of $\C^{nN}$ having dimension $\max\{r \mid D_{\g, r}(\a^{-1})=0\}$. 

Assume that $\a$ is a zero of an irreducible factor of $\D_\g(t)$ that does not occur in $\D_{\g, 2}(t)$. Since ${\rm Hom}(\pi_\g, {\cal D}_N)_\a^0$ has dimension 1, there exists an index $i$ such that projection onto some coordinate $b_i^{(j)}$ of $b_i\in \C^N$ is an isomorphism. One checks that ${\rm Hom}(\pi_\g, {\cal D}_N)_\a^0$ is invariant under conjugation by $d: z \mapsto \beta z$, for any $\beta$. Moreover, $b_i^{(j)}$ is transformed into $\beta^{-1} b_i^{(j)}$. As $\beta$ ranges over all nonzero complex numbers, so does $b_i^{(j)}$. \end{proof}

In view the theorem of Burde and de Rham, it is natural to ask
whether  the set of zeros of the Alexander-Lin polynomial $D_\g(t)$ is necessarily closed under inversion. The question was originally posed in \cite{kitano} for the $\D_\g(t)$. It is relatively easy to construct examples that provide a negative answer both for $D_\g(t)$ and $\D_\g(t)$ by mapping meridians of  the knot group $\pi$ to matrices with determinant other than $\pm 1$. (See \cite{hsw} for a particularly simple example suggested by S. Friedl.) An example of a representation $\g: \pi \to {\rm SL}_3\Z$ for which $D_\g(t)$ is non-reciprocal is given in \cite{hsw}.


 \bigskip

\end{document}